\documentclass[reqno,11pt]{amsart}
\usepackage[foot]{amsaddr}
\usepackage{amssymb,amsmath,amsthm,amsfonts}
\usepackage{mathrsfs,dsfont,comment,mathscinet,mathtools}
\usepackage{regexpatch}
\usepackage{graphicx,pgfplots}
\usepackage{enumitem}

\usepackage[utf8]{inputenc}
\usepackage[normalem]{ulem}
\usepackage[left=2.5cm,right=2.5cm,top=2.5cm,bottom=2.5cm]{geometry}

\usepackage{pifont}
\usepackage{tkz-fct}

\usepackage{graphicx,tikz,color}
\usetikzlibrary{quotes,angles}

\usepackage[format=hang,labelfont=bf]{caption}

\usepackage{esint}

\usepackage{ wasysym }

\usepackage{epstopdf}
\usepackage[colorlinks=true, pdfstartview=FitV, linkcolor=blue, 
            citecolor=blue, urlcolor=blue]{hyperref}
            
\allowdisplaybreaks
\numberwithin{equation}{section}
\theoremstyle{plain}
\newtheorem{theorem}{Theorem}[section]
\newtheorem{proposition}[theorem]{Proposition}
\newtheorem{lemma}[theorem]{Lemma}

\newtheorem{definition}[theorem]{Definition}
\theoremstyle{definition}
\newtheorem{remark}[theorem]{Remark}

\newtheorem*{theorem*}{Theorem}

\definecolor{mblue}{HTML}{13439b}

\newcommand{\N}{\mathbb{N}}
\newcommand{\R}{\mathbb{R}}
\newcommand{\RN}{\mathbb{R}^N}
\newcommand{\RNN}{\mathbb{R}^{N\times N}}

\newcommand{\SN}{\mathbb{S}^{N-1}}

\newcommand{\W}{\mathrm{W}}

\renewcommand{\d}{\mathrm{d}}
\newcommand{\restr}[1]{|_{#1}}

\def\Xint#1{\mathchoice
{\XXint\displaystyle\textstyle{#1}}%
{\XXint\textstyle\scriptstyle{#1}}%
{\XXint\scriptstyle\scriptscriptstyle{#1}}%
{\XXint\scriptscriptstyle\scriptscriptstyle{#1}}%
\!\int}
\def\XXint#1#2#3{{\setbox0=\hbox{$#1{#2#3}{\int}$ }
\vcenter{\hbox{$#2#3$ }}\kern-.6\wd0}}

\def\dashint{\Xint-}

\makeatletter
\newcommand{\subalign}[1]{%
	\vcenter{%
		\Let@ \restore@math@cr \default@tag
		\baselineskip\fontdimen10 \scriptfont\tw@
		\advance\baselineskip\fontdimen12 \scriptfont\tw@
		\lineskip\thr@@\fontdimen8 \scriptfont\thr@@
		\lineskiplimit\lineskip
		\ialign{\hfil$\m@th\scriptstyle##$&$\m@th\scriptstyle{}##$\hfil\crcr
			#1\crcr
		}%
	}%
}
\makeatother

\newcommand{\mres}{\mathbin{\vrule height 1.6ex depth 0pt width
		0.13ex\vrule height 0.13ex depth 0pt width 1.3ex}}

\makeatletter
\newcommand{\subsetcong}{\mathrel{\mathpalette\subseteq@cong\relax}}
\newcommand{\subsetsim}{\mathrel{\mathpalette\subset@sim\relax}}

\newcommand{\subseteq@cong}[2]{%
	\vbox{\offinterlineskip\m@th
		\ialign{\hfil$#1##$\hfil\cr
			\sim\cr\subseteq\cr
		}%
	}%
}

\newcommand{\subset@sim}[2]{%
	\vbox{\offinterlineskip\m@th
		\ialign{\hfil$#1##$\hfil\cr
			\sim\cr\subset\cr
		}%
	}%
}

\makeatother

\newcommand{\tr}{\mathrm{tr}}
\newcommand{\cof}{\mathrm{cof}}
\newcommand{\adj}{\mathrm{adj}}

\newcommand{\per}{\mathrm{Per}}
\renewcommand{\div}{\mathrm{div}}
\newcommand{\dist}{\mathrm{dist}}
\newcommand\wk{\rightharpoonup}
\newcommand{\wks}{\overset{\ast}{\rightharpoonup}}

\newcommand*\closure[1]{\overline{#1}}

\newcommand{\leb}{\mathscr{L}^N}
\newcommand{\haus}{\mathscr{H}^{N-1}}

\newcommand{\imt}{\mathrm{im}_{\rm T}}
\newcommand{\img}{\mathrm{im}_{\rm G}}
\newcommand{\domg}{\mathrm{dom}_{\rm G}}

\definecolor{cardinal}{rgb}{0.77, 0.12, 0.23}

\newcommand\EEE{\color{black}}

\newcommand{\MMM}{\color{black}}

\makeatletter
\renewcommand{\subsetsim}{\mathrel{\mathpalette\subset@cong\relax}}

\newcommand{\subset@cong}[2]{%
	\vbox{\offinterlineskip\m@th
		\ialign{\hfil$#1##$\hfil\cr
			\sim\cr\subset\cr
		}%
	}%
}

\makeatletter
\renewcommand{\subsetcong}{\mathrel{\mathpalette\subset@sim\relax}}

\renewcommand{\subset@sim}[2]{%
	\vbox{\offinterlineskip\m@th
		\ialign{\hfil$#1##$\hfil\cr
			\sim\cr\subseteq\cr
		}%
	}%
}

\makeatletter
\xpatchcmd{\proof}{\itshape}{\bfseries}{}{}
\makeatother

\title[Anisotropic cavitation energies]{Anisotropic energies for the \MMM modeling \EEE of cavitation in nonlinear elasticity}
\author[M. Bresciani]{Marco Bresciani${}^{*}$}
\address{* Department of Mathematics, Friedrich-Alexander-Universit\"{a}t Erlangen-N\"{u}rnberg, Cauerstrasse 11, 91058 Erlangen, Germany}
\email{marco.bresciani@fau.de}

\date{\today}

\begin{document}

\setlength\parindent{0pt}

\vskip .2truecm
\begin{abstract}
   We {\MMM study} free-discontinuity functionals in nonlinear elasticity, where discontinuities correspond to the phenomenon of cavitation. The energy comprises  two terms: a volume term accounting for the elastic energy;   and a surface term  concentrated on the boundaries of the cavities in the deformed configuration that depends on their  unit normal. {\MMM First, we prove the existence of energy-minimizing deformations.} While the treatment of the {\MMM volume} term is standard,  that of the {\MMM surface} term relies on the regularity of inverse deformations, their weak continuity properties, and {\MMM Ambrosio's} lower semicontinuity theorem for special functions of bounded variation. {\MMM Additionally, we identify  sufficient conditions for minimality by employing outer variations and applying  the formula for the first derivative of the anisotropic perimeter.}
\end{abstract}
\maketitle

\section{Introduction}

\subsection{Motivation and state of the art} Modern theories of mechanics (see, e.g., \cite{cho.gent,gent.wang,lopezpamies1,william.sharpery}) suggest that cavitation in  solids originates from the competition between two contributions: a volume energy accounting for the elastic response of materials;  and surface energy governing the growth of  cavities. Here, following {\sc Ball} \cite{ball.op}, we refer to cavitation as the  sudden formation of voids inside of elastic materials in response to tensile stresses. 

A first variational theory of cavitation in this spirit  has been advanced by {\sc M\"{u}ller \& Spector}. In their seminal paper \cite{mueller.spector}, they investigated  existence and  regularity of minimizers of the functional
\begin{equation}
	\label{eq:intro-MS}
	\boldsymbol{y}\mapsto \int_\Omega \MMM W( D\boldsymbol{y})\EEE \,\d\boldsymbol{x}+\per \left( \boldsymbol{y}(\Omega) \right),
\end{equation} 
among all deformations $\boldsymbol{y}\colon \Omega\to \RN$, with $\Omega\subset \RN$ denoting the reference configuration of an elastic object, belonging to the Sobolev space $W^{1,p}(\Omega;\RN)$ with $p>N-1$.
The first term in \eqref{eq:intro-MS} involves a material energy density $W$ satisfying standard assumptions of nonlinear elasticity, while the second one measures the perimeter of the deformed configuration (to be interpreted in a precise measure-theoretic sense). We mention  \cite{conti.delellis} for a related study with deformation in $W^{1,N-1}(\Omega;\RN)$. 

Subsequently, {\sc Henao \& Mora-Corral} \cite{henao.moracorral.invertibility} have addressed the minimization of the functional
\begin{equation}
	\label{eq:intro-HMC1}
	 \boldsymbol{y}\mapsto \int_\Omega \MMM W(D\boldsymbol{y})\EEE \,\d\boldsymbol{x} + \mathcal{S}(\boldsymbol{y})
\end{equation}
for deformations $\boldsymbol{y}\in W^{1,N-1}(\Omega;\RN)$. 
As before, the first term in \eqref{eq:intro-HMC1} accounts for the elastic energy. The expression of the second term  is quite involved (see Definition \ref{def:surface-energy} below) but finds a complete justification within the theory of cartesian currents \cite{gms.cc1,gms.cc2}. Roughly speaking,  this term measures the new surface created by the deformation in the sense that
\begin{equation}
	\label{eq:intro-S}
	\mathcal{S}(\boldsymbol{y})=\per \left( \boldsymbol{y}(\Omega) \right)-\haus(\boldsymbol{y}(\partial \Omega)),
\end{equation} 
where $\haus$ denotes the $(N-1)$-dimensional Hausdorff measure in $\RN$. The interpretation of $\mathcal{S}(\boldsymbol{y})$ as the area of the new surface created by $\boldsymbol{y}$ was later made rigorous in terms of geometric measure theory by {\sc Henao \& Mora-Corral}   \cite{henao.moracorral.fracture}. 

Although the study in \cite{henao.moracorral.invertibility} was mainly motivated by the search of a unified variational framework for cavitation and fracture, later on 
  {\sc Henao \& Mora-Corral} \cite{henao.moracorral.lusin} have shown that when fracture is excluded, i.e., in the case of Sobolev deformations (see Remark \ref{rem:inv-finite} below),  the energy in \eqref{eq:intro-HMC1} can be rewritten as
\begin{equation}
	\label{eq:intro-HMC2}
	 \boldsymbol{y}\mapsto \int_\Omega \MMM  W(D\boldsymbol{y})\EEE \,\d\boldsymbol{x} + \sum_{\boldsymbol{a}\in C_{\boldsymbol{y}}} \per \left( \imt(\boldsymbol{y},\boldsymbol{a})  \right).
\end{equation}
In the previous equation, $C_{\boldsymbol{y}}$ stands for the set of points in correspondence of which the deformation opens a cavity with $\imt(\boldsymbol{y},\boldsymbol{a})$ being the cavity opened at $\boldsymbol{a}\in\Omega$ and $\per \left(\imt(\boldsymbol{y},\boldsymbol{a})\right)$ its perimeter. The formulation in \eqref{eq:intro-HMC2} appears evidently as  a free-discontinuity problem given that the set $C_{\boldsymbol{y}}$ is not prescribed (see   \cite{henao,sivaloganathan.spector,sivaloganathan.spector.tilakraj} for some work in that direction) but rather selected by minimality. Bearing \eqref{eq:intro-S} in mind, the energy in \eqref{eq:intro-HMC2} appears more physical compared with that in \eqref{eq:intro-MS} as it does not penalize the stretching of the outer boundary. An  evolutionary version of the model in \cite{henao.moracorral.lusin} has been studied by {\sc Mora-Corral} \cite{moracorral} in the quasistatic setting.

\MMM Concerning related developments, disparate approaches for the numerical approximation of the functional in \eqref{eq:intro-HMC2} have been investigated by means of $\Gamma$-convergence in \cite{bresciani.friedrich.core,henao.moracorral.xu}. Eventually, we mention an extension of the model in \cite{henao.moracorral.invertibility,henao.moracorral.lusin} for active materials (e.g.,  nematic and magnetic elastomers)  investigated in \cite{bresciani.friedrich.moracorral}.
 \EEE

\subsection{Contributions of the paper}
In this work, we propose  a generalization of the model in \cite{henao.moracorral.lusin} accounting for the effect of material anisotropy or, said differently, of the geometry of the cavities. Specifically, we \MMM are interested in \EEE minimizers of the  functional
\begin{equation}
	\label{eq:intro-B}
	\boldsymbol{y}\mapsto \int_\Omega \MMM  W(D\boldsymbol{y})\EEE \,\d \boldsymbol{x}+\sum_{\boldsymbol{a}\in C_{\boldsymbol{y}}} \int_{\partial^* \imt(\boldsymbol{y},\boldsymbol{a})} \phi(\boldsymbol{\nu}_{\imt(\boldsymbol{y},\boldsymbol{a})})\,\d\haus
\end{equation}
among deformations $\boldsymbol{y}\in W^{1,p}(\Omega;\RN)$. The choice of the function space is merely technical and it is mainly motivated by the continuity of the topological degree (recall that $p>N-1$), which does not hold in $W^{1,N-1}(\Omega;\RN)$.  Also here, the first term in \eqref{eq:intro-B} stands for the elastic energy. The second term is given by the integral of a density $\phi$ on the (reduced) boundary of the cavities with $\phi$ depending on the   unit normal $\boldsymbol{\nu}_{\imt(\boldsymbol{y},\boldsymbol{a})}$ to each cavity $\imt(\boldsymbol{y},\boldsymbol{a})$ for $\boldsymbol{a}\in C_{\boldsymbol{y}}$. This surface {\MMM energy} is reminiscent of free-discontinuity functionals in the space of special functions of bounded variations, where the  energy generally depends on the normal to the jump set \cite{ambrosio,ambrosio.braides}. In particular, taking {\MMM $\phi(\boldsymbol{z})=|\boldsymbol{z}|$} in \eqref{eq:intro-B},  we recover \eqref{eq:intro-HMC2}.

\MMM Our main results  are Theorem \ref{thm:existence} and Theorem \ref{thm:condition} below. The first theorem establishes the existence of minimizers for the energy \eqref{eq:intro-B} among admissible deformations satisfying Dirichlet boundary conditions, 
 while the second one identifies some sufficient conditions for minimality. These findings provide an extension of analogous results in \cite{henao.moracorral.fracture,henao.moracorral.lusin,mueller.spector} to the anisotropic setting for deformations in $W^{1,p}(\Omega;\RN)$.

For the proof of Theorem \ref{thm:existence}, 
 \EEE 
the key idea for treating the surface term in \eqref{eq:intro-B} is to rewrite it in terms of a suitably defined inverse deformation $\widehat{\boldsymbol{y}}^{-1}$  of ${\boldsymbol{y}}$ (see Definition \ref{def:inverse-sobolev} below).  Taking advantage of a regularity result by {\sc Henao \& Mora-Corral} \cite{henao.moracorral.regularity}, we infer that $\widehat{\boldsymbol{y}}^{-1}$ is a special function of bounded variation on the topological image of $\boldsymbol{y}$ (see Definition \ref{def:topim} below). Then, the convergence properties of inverse deformations established in \cite{bresciani.friedrich.moracorral} and a classical lower semicontinuity theorem due to {\sc Ambrosio} \cite{ambrosio} allow for a successful application of the direct method of the calculus of variations. \MMM  The proof of Theorem \ref{thm:condition} combines the classical approach based on outer variations \cite{ball.op} with an anisotropic version of the first-variation formula for the perimeter \cite{bellettini.novaga.riey}. \EEE

\MMM In conclusion, \EEE we observe that the functional \eqref{eq:intro-B} depends only on the first-order geometry of the cavities. It would be very interesting to consider a possible dependence also on their second-order geometry, i.e., on their curvature. However, even a proper formulation of a such problem  is currently  unclear.  

 \subsection{Structure of the paper}  The paper is organized into four sections, the first one being this Introduction. In Section \ref{sec:main-result}, we specify our setting and we state our main results. \EEE  In Section \ref{sec:prelim}, we develop all the necessary tools for the \MMM proofs of our main results \EEE which \MMM are \EEE  presented in Section \ref{sec:proof}.

\subsection*{Notation} We employ standard notation for measure theory and  function spaces. We write $B(\boldsymbol{a},r)$ and $S(\boldsymbol{a},r)$ for the open ball and the sphere centered at $\boldsymbol{a}\in\RN$ with radius $r>0$. 
The Lebesgue measure and the $(N-1)$-dimensional Hausdorff measure on $\RN$ are denoted by $\leb$ and $\haus$, respectively. Unless differently stated, the expression ``almost everywhere'' (abbreviated as ``a.e.'') and all similar ones are  referred to $\leb$. Given two sets $E,F\subset \RN$, the notation $E\subset \subset F$ means that $\closure{E}\subset F$, where $\closure{E}$ is the closure of $E$. Also, we write $E\cong F$ and $E\simeq F$ whenever $\leb(E \triangle F)=0$ and $\haus(E \triangle F)=0$, respectively. Here, $E\triangle F\coloneqq (E\setminus F)\cup (F\setminus E)$. Given two measurable functions $\boldsymbol{u}$ and $\boldsymbol{v}$, we write $\boldsymbol{u}\cong \boldsymbol{v}$ for $\leb(\{\boldsymbol{u}\neq \boldsymbol{v}\})=0$ and $\boldsymbol{u}\simeq \boldsymbol{v}$ for $\haus(\{\boldsymbol{u}\neq \boldsymbol{v}\})=0$. We employ standard notation for sets of finite perimeter and for functions of bounded variation. In particular, for a set $E\subset \RN$ with  finite perimeter, we denote by $\partial^*E$ and $\per(E)$ its reduced boundary and perimeter, respectively, and, for a function $\boldsymbol{u}$ of bounded variation, we indicate by $J_{\boldsymbol{u}}$ its set of jump points {\MMM and by $\boldsymbol{\nu}_{\boldsymbol{u}}$ the associated unit normal. }

\section{Setting and main result}
\label{sec:main-result}

Let $N\in\N$ with $N\geq 2$ and $\Omega\subset \RN$ be a bounded Lipschitz domain. Fix $p>N-1$ and consider the class of admissible deformations
\begin{equation*}
	\mathcal{Y}_p(\Omega)\coloneqq \left\{ \boldsymbol{y}\in W^{1,p}(\Omega;\RN):\; \text{$\det D \boldsymbol{y}(\boldsymbol{x})>0$ for a.e. $\boldsymbol{x}\in\Omega$},\: \text{$\boldsymbol{y}$ satisfies (INV)} \right\}.
\end{equation*}
In the previous equation, we refer to condition (INV)  by {\sc M\"{u}ller \& Spector} \cite{mueller.spector} recalled in Definition \ref{def:INV} below. Roughly speaking, this condition prohibits that cavities created at one point can be filled with material coming from elsewhere. 

We are interested in the minimization of the following energy
\begin{equation}
	\label{eq:E}
	\mathcal{E}(\boldsymbol{y})\coloneqq \int_\Omega \MMM  W(D\boldsymbol{y})\EEE \,\d \boldsymbol{x}+\sum_{\boldsymbol{a}\in C_{\boldsymbol{y}}} \int_{\partial^* \imt(\boldsymbol{y},\boldsymbol{a})} \phi(\boldsymbol{\nu}_{\imt(\boldsymbol{y},\boldsymbol{a})})\,\d\haus.
\end{equation} 
The first term in \eqref{eq:E} features the bulk density \MMM $W\colon  \RNN_+ \to \R$, where $\RNN_+$ denotes the set of matrices $\boldsymbol{F}\in\RNN$ with $\det \boldsymbol{F}>0$, \EEE on which we assume the following:
\begin{enumerate}[label=(W\arabic*)]
	\item \emph{Coercivity:} \MMM There exist a constant $c>0$ and \EEE  a Borel function $\gamma \colon (0,+\infty)\to [0,+\infty)$ with
	\begin{equation}
		\label{eq:gamma}
		\lim_{h\to 0^+}\gamma(h)=+\infty, \quad \lim_{h\to +\infty} \frac{\gamma(h)}{h}=+\infty
	\end{equation}
	such that
	\begin{equation}
		\label{eq:growth}
		\text{\MMM $W(\boldsymbol{F})\geq c|\boldsymbol{F}|^p+\gamma(\det \boldsymbol{F})$ \quad for all $\boldsymbol{F}\in\RNN_+$;}
	\end{equation}
	\item \emph{Polyconvexity:} There exists a \MMM convex function $g\colon  \prod_{i=1}^{N-1} \R^{\binom{N}{i}\times \binom{N}{i}} \times (0,+\infty) \to \R$  
	 such that
	\begin{equation*}
		W(\boldsymbol{F})=g(\mathbf{M}(\boldsymbol{F})) \quad \text{for all $\boldsymbol{F}\in\RNN_+$,}
	\end{equation*}
	where $\mathbf{M}(\boldsymbol{F})$ denotes  the collection of all minors of $\boldsymbol{F}$ arranged in a given order;
	\item \emph{Regularity and control:} The function $W$ is continuously differentiable and there exists a constant $\tilde{c}>0$ such that
	\begin{equation*}
		|DW(\boldsymbol{F})\boldsymbol{F}^\top|\leq \tilde{c} (W(\boldsymbol{F})+1) \quad \text{for all $\boldsymbol{F}\in\RNN_+$.}
	\end{equation*}
\end{enumerate}
\MMM These assumptions are standard in nonlinear elasticity.  Generally, (W1)--(W2) are assumed for existence results, while (W3) is employed to derive first-order conditions. Note that $W$ takes only finite values, so that (W1)--(W2) yield the continuity of $W$.  \EEE

The second term in \eqref{eq:E} is given by a sum over the set $C_{\boldsymbol{y}}$ of cavitation points of $\boldsymbol{y}$. This set is introduced in Definition \ref{def:cav} below. For each $\boldsymbol{a}\in C_{\boldsymbol{y}}$, the set $\imt(\boldsymbol{y},\boldsymbol{a})$ is the topological image of $\boldsymbol{a}$ under $\boldsymbol{y}$ representing the cavity opened by $\boldsymbol{y}$ in correspondence of the point $\boldsymbol{a}$. 
On the surface density \MMM $\phi\colon \RN \to (0,+\infty)$ \EEE we make the following  assumptions:
\begin{enumerate}[label=($\Phi$\arabic*)]
	\item \MMM \emph{Homogeneity:} $\phi$ is one-homogeneous, namely,
	\begin{equation*}
		\phi(\boldsymbol{z})=|\boldsymbol{z}| \phi \left( \frac{\boldsymbol{z}}{|\boldsymbol{z}|} \right) \quad \text{for all $\boldsymbol{z}\in\RN$;}
	\end{equation*}
	\item \emph{Lower bound:} There exists a constant $\hat{c}>0$ such that 
	\begin{equation*}
		\text{$\phi(\boldsymbol{z})\geq \hat{c}|\boldsymbol{z}|$ \quad for all $\boldsymbol{z}\in\RN$; }
	\end{equation*} \EEE 
	\item \emph{Convexity:} $\phi$ is convex; 
	\item \MMM \emph{Regularity:} $\phi$ is continuously differentiable in a neighborhood of $\SN$, where $\SN\subset \RN$ denotes the unit sphere centered at the origin.
\end{enumerate}
 Assumptions ($\Phi$1)--($\Phi$3) are standard in the existence theory of free-discontinuity problems for special functions of bounded variation. Precisely, these assumptions characterize biconvex integrands solely depending  on the normal to the jump set \cite{ambrosio,ambrosio.braides}. \MMM Note that ($\Phi$3) ensures the continuity of $\phi$. The regularity in ($\Phi$4) will be only relevant while deriving sufficient conditions for minimality.   \EEE 

\MMM We introduce the following class of admissible deformations
\begin{equation}
	\label{eq:A}
	\mathcal{A}\coloneqq\left\{ \boldsymbol{y}\in\mathcal{Y}_p(\Omega):\:\mathcal{S}(\boldsymbol{y})<+\infty,\:\text{$\boldsymbol{y}=\boldsymbol{d}$ on $\Gamma$ in the sense of traces}  \right\}.
\end{equation} 
In this equation, the boundary datum $\boldsymbol{d}\colon \Gamma\to \RN$ and the set  $\Gamma\subset \partial \Omega$  with $\haus(\Gamma)>0$ are both $\haus$-measurable. The functional $\mathcal{S}$, originally introduced in \cite{henao.moracorral.invertibility}, is given in Definition \ref{def:surface-energy} below. In this regard,  we mention that the condition $\mathcal{S}(\boldsymbol{y})<+\infty$  constitutes a very mild regularity requirement which is violated only by deformations exhibiting pathological behaviors.  \EEE 

\MMM Our first main result asserts the existence of minimizers among admissible deformations. \EEE 

\begin{theorem}[Existence of minimizers]
	\label{thm:existence}
Assume {\rm (W1)--(W2)} and {\rm ($\Phi$1)--($\Phi$3)} and  suppose that the class $\mathcal{A}$ in \eqref{eq:A} is nonempty	.
Then, the functional $\mathcal{E}$ in \eqref{eq:E} admits minimizers among all deformations in   $\mathcal{A}$.
\end{theorem}

\MMM The previous theorem admits variants for inhomogeneous bulk densities as well as for incompressible materials. In both cases, the adaptation of the proof is straightforward.

Our second main result provides sufficient conditions for the minimality of admissible deformations and reads as follows.

\begin{theorem}[Sufficient conditions for minimality]
	\label{thm:condition}
	Assume {\rm (W3)}, {\rm ($\Phi$2)}, and {\rm ($\Phi$4)}. Suppose that $\boldsymbol{y}$ is a minimizers of the functional $\mathcal{E}$ in \eqref{eq:E} among all deformations of the class  $\mathcal{A}$ in \eqref{eq:A} with $\mathcal{E}(\boldsymbol{y})<+\infty$. Then, for any $\boldsymbol{\psi}\in C^1_{\rm c}(\RN;\RN)$ such that $\boldsymbol{\psi}\circ \boldsymbol{d}=\boldsymbol{0}$ on $\Gamma$ in the sense of traces, there holds
	\begin{equation*}
		\label{eq:condition}
		\begin{split}
			\int_\Omega &\left( DW(D\boldsymbol{y})(D\boldsymbol{y})^\top  \right):D\boldsymbol{\psi}(\boldsymbol{y})\,\d\boldsymbol{x}
			+\sum_{\boldsymbol{a}\in C_{\boldsymbol{y}}} \int_{\partial^*\imt(\boldsymbol{y},\boldsymbol{a})} \phi(\boldsymbol{\nu}_{\imt(\boldsymbol{y}_n,\boldsymbol{a})})\,\div_{\phi}^{\partial^*\imt(\boldsymbol{y},\boldsymbol{a})}\boldsymbol{\psi}\,\d\haus=0,
		\end{split}
	\end{equation*}
	where 
	\begin{equation*}
		\div_{\phi}^{\partial^*\imt(\boldsymbol{y},\boldsymbol{a})} \boldsymbol{\psi}\coloneqq \div\, \boldsymbol{\psi}-D\phi(\boldsymbol{\nu}_{\imt(\boldsymbol{y},\boldsymbol{a})}) \cdot \frac{(D\boldsymbol{\psi})^\top \boldsymbol{\nu}_{\imt(\boldsymbol{y},\boldsymbol{a})}   }{\phi(\boldsymbol{\nu}_{\imt(\boldsymbol{y},\boldsymbol{a})})}. 
	\end{equation*}
\end{theorem}
When $\phi$ is twice continuously differentiable in a neighborhood of $\SN$, $\boldsymbol{y}$ creates a finite number of cavities each of which has regular boundary, and $\boldsymbol{y}$ is a smooth diffeormorphism on $\closure{\Omega}\setminus C_{\boldsymbol{y}}$, an anisotropic version of the Gauss theorem  on surfaces (see, e.g., \cite[Equation (22)]{bellettini.novaga.paolini}) allows us to recognize
the conditions  in Theorem \ref{thm:condition} as the weak formulation of the boundary-value problem 
\begin{equation*}
	\begin{cases}
		\div\,\boldsymbol{T}^{\boldsymbol{y}}=\boldsymbol{0} & \text{in $\boldsymbol{y}(\Omega \setminus C_{\boldsymbol{y}})$,}\\
		\boldsymbol{T}^{\boldsymbol{y}}\boldsymbol{\nu}_{\imt(\boldsymbol{y},\boldsymbol{a})}=-h^\phi_{\imt(\boldsymbol{y},\boldsymbol{a})}\boldsymbol{\nu}_{\imt(\boldsymbol{y},\boldsymbol{a})} & \text{on $\partial\, \imt(\boldsymbol{y},\boldsymbol{a})$ for all $\boldsymbol{a}\in C_{\boldsymbol{y}}$,}\\
	\end{cases}
\end{equation*} 
where  $\boldsymbol{T}^{\boldsymbol{y}}\colon \boldsymbol{y}(\Omega \setminus C_{\boldsymbol{y}})\to \RNN$ and  $h^\phi_{\imt(\boldsymbol{y},\boldsymbol{a})}\colon \partial\,\imt(\boldsymbol{y},\boldsymbol{a}) \to \R$ are the Cauchy stress tensor and the anisotropic mean curvature, respectively, defined as $\boldsymbol{T}^{\boldsymbol{y}}(\boldsymbol{y})\coloneqq \left(DW(D\boldsymbol{y})(D\boldsymbol{y})^\top  \right)(\det D \boldsymbol{y})^{-1}$ and $h^\phi_{\imt(\boldsymbol{y},\boldsymbol{a})}\coloneqq \tr \left( D^2\phi(\boldsymbol{\nu}_{\imt(\boldsymbol{y},\boldsymbol{a})}) D\boldsymbol{\nu}_{\imt(\boldsymbol{y},\boldsymbol{a})}  \right)$.

To establish Theorem \ref{thm:condition}, we employ the classical method based on outer variations. Considering inner variations is also possible. However, as topological images are invariant under this kind of variations, this approach does not provide any information on the surface energy, so that one simply recovers well-known sufficient conditions for the elastic energy (see, e.g.,  \cite[Theorem 2.3(ii)]{ball.op} and  \cite[Theorem 6.2]{mueller.spector}). 

\EEE  

\section{Continuity properties and regularity of inverse deformations}
\label{sec:prelim}

In this section, we collects some results that will be instrumental for the proof of Theorem \ref{thm:existence}. Having future applications in mind, we distinguish the two cases of approximately differentiable and Sobolev deformations. 

\subsection{Approximately differentiable deformations}
For the notion of approximate differentiability, we refer to \cite[Section 3.1]{gms.cc1}.  To ease the exposition, we introduce the following class of deformations
\begin{equation*}
	\begin{split}
		\mathcal{Y}(\Omega)\coloneqq \big \{ \boldsymbol{y}\colon \Omega \to \RN:\hspace{4pt}  \text{$\boldsymbol{y}$ a.e. approximately differentiable},\:\text{$\boldsymbol{y}$ a.e. injective,}& \\ \text{$\det \nabla \boldsymbol{y}(\boldsymbol{x})>0$ for a.e. $\boldsymbol{x}\in\Omega$}& \big \}.
	\end{split}
\end{equation*}
Henceforth, we denote by $\nabla \boldsymbol{y}$  the approximate gradient of $\boldsymbol{y}$.
The assumption $\boldsymbol{y}$ almost everywhere injective requires the existence of a set $X\subset \Omega$ with $\leb(X)=0$ such that  $\boldsymbol{y}\restr{\Omega \setminus X}$ is injective.

The following functional, measuring the area of the new surface created by deformations, has been introduced by {\sc Henao \& Mora-Corral} \cite{henao.moracorral.invertibility} in their study on the weak continuity of Jacobian determinants. 

\begin{definition}[Surface energy functional]
	\label{def:surface-energy}
	Given $\boldsymbol{y}\in\mathcal{Y}(\Omega)$,	we define the funcitonal
	\begin{equation*}
		\mathcal{S}(\boldsymbol{y})\coloneqq \left\{ \mathcal{S}_{\boldsymbol{y}}(\boldsymbol{\eta}):\:\boldsymbol{\eta}\in C^0_{\rm c}(\Omega\times \RN;\RN),\:\|\boldsymbol{\eta}\|_{L^\infty(\Omega\times \RN;\RN)}\leq 1   \right\},
	\end{equation*}
	where
	\begin{equation*}
		\mathcal{S}_{\boldsymbol{y}}(\boldsymbol{\eta})\coloneqq \int_\Omega \left\{(\cof D \boldsymbol{y}(\boldsymbol{x}))\cdot D_{\boldsymbol{x}}\boldsymbol{\eta}(\boldsymbol{x},\boldsymbol{y}(\boldsymbol{x}))+\div_{\boldsymbol{\xi}}\boldsymbol{\eta}(\boldsymbol{x},\boldsymbol{y}(\boldsymbol{x}))\det D \boldsymbol{y}(\boldsymbol{x}) \right\}\,\d\boldsymbol{x}
	\end{equation*}
	with $(\boldsymbol{x},\boldsymbol{\xi})\in\Omega\times \RN$ denoting the variables of $\boldsymbol{\eta}$.
\end{definition}

We introduce the notions of geometric domain and image as in \cite[Definition 2.4]{bresciani.friedrich.moracorral}. For the notion of density of a set with respect to a point, we refer to \cite[Section 3.1]{gms.cc1}

\begin{definition}[Geometric domain and image]
	Given $\boldsymbol{y}\in\mathcal{Y}(\Omega)$, we define the  geometric domain of $\boldsymbol{y}$ as the set $\domg(\boldsymbol{y},\Omega)$   of points $\boldsymbol{x}\in \Omega$ such that  $\boldsymbol{y}$ is approximately differentiable at $\boldsymbol{x}$  with $\det\nabla \boldsymbol{y}(\boldsymbol{x})>0$, and  there exist a compact set $K\subset \Omega$ having density one at $\boldsymbol{x}$ and function $\boldsymbol{w}\in C^1(\RN;\RN)$ such that $\boldsymbol{w}\restr{K}=\boldsymbol{y}\restr{K}$ and $(\nabla \boldsymbol{y})\restr{K}=(D\boldsymbol{w})\restr{K}$.
	Moreover, we define the geometric image of $\boldsymbol{y}$ as the set $\img(\boldsymbol{y},\Omega)\coloneqq \boldsymbol{y}(\domg(\boldsymbol{y},\Omega))$.
\end{definition}

Following \cite[Proposition 2.4]{mueller.spector}, it is shown in \cite[p. 585]{henao.moracorral.fracture} that the set $\domg(\boldsymbol{y},\Omega)$ has full measure in $\Omega$.  Also, in \cite[Lemma 3]{henao.moracorral.fracture}, it is proved that the restriction $\boldsymbol{y}\restr{\domg(\boldsymbol{y},\Omega)}\colon \domg(\boldsymbol{y},\Omega)\to \RN$ is (everywhere) injective. 
This fact allows us to  rigorously define inverse deformations. 

\begin{definition}[Inverse deformation]
	\label{def:inverse}
Given $\boldsymbol{y}\in\mathcal{Y}(\Omega)$, we define the inverse deformation $\boldsymbol{y}^{-1}\colon \img(\boldsymbol{y},\Omega)\to \RN$ of $\boldsymbol{y}$ by setting $\boldsymbol{y}^{-1}(\boldsymbol{\xi})\coloneq \boldsymbol{x}$, where $\boldsymbol{x}$ is the unique point in $\domg(\boldsymbol{y},\Omega)$ such that $\boldsymbol{y}(\boldsymbol{x})=\boldsymbol{\xi}$. Moreover, we define its extension $\overline{\boldsymbol{y}}^{-1} \colon \RN \to \RN$ by setting
\begin{equation*}
	\overline{\boldsymbol{y}}^{-1}(\boldsymbol{\xi})\coloneqq \begin{cases}
		\boldsymbol{y}^{-1}(\boldsymbol{\xi}) & \text{if $\boldsymbol{\xi}\in \img(\boldsymbol{y},\Omega)$,}\\
		\boldsymbol{0} & \text{if $\boldsymbol{\xi}\in \RN\setminus \img(\boldsymbol{y},\Omega)$.}
	\end{cases} 
\end{equation*} 
\end{definition}

By \cite[Proposition 1(ii), p.221]{gms.cc1}, the restriction $\boldsymbol{y}\restr{\domg(\boldsymbol{y},\Omega)}$ maps measurable sets to measurable sets. Thus,  $\boldsymbol{y}^{-1}$ is a measurable map. Additionally, 
as shown in \cite[Lemma 3]{henao.moracorral.xu}, the inverse $\boldsymbol{y}^{-1}$ is (everywhere) approximately differentiable with $\nabla \boldsymbol{y}^{-1}=(\nabla \boldsymbol{y})^{-1}\circ \boldsymbol{y}^{-1}$ in $\img(\boldsymbol{y},\Omega)$, while for its extension we have $\nabla \overline{\boldsymbol{y}}^{-1}=\nabla \boldsymbol{y}^{-1}$ in $\img(\boldsymbol{y},\Omega)$ and $\nabla \overline{\boldsymbol{y}}^{-1}\cong \boldsymbol{O}$ in $\RN\setminus \img(\boldsymbol{y},\Omega)$.

The continuity properties of the inverse for approximately differentiable deformations have been examined in \cite[Proposition 3.3]{bresciani.friedrich.moracorral}. Specifically, the weak convergence of approximate gradients of inverse deformations is stated without proof in \cite[Remark 3.3]{bresciani.friedrich.moracorral}.  

\begin{lemma}[Continuity properties of inverse deformations]
	\label{lem:continuity-inverse}
	Let $(\boldsymbol{y}_n)$ be a sequence in $\mathcal{Y}(\Omega)$ and $\boldsymbol{y}\in\mathcal{Y}(\Omega)$. Suppose that 
	\begin{equation}
		\label{eq:conv-ae-det}
		\text{$\boldsymbol{y}_n \to \boldsymbol{y}$ a.e. in $\Omega$,}\quad \text{$\det \nabla \boldsymbol{y}_n \wk \det \nabla \boldsymbol{y}$ in $L^1(\Omega)$.}
	\end{equation}
	Then
	\begin{align}
		\label{eq:conv-img}
		\chi_{\img(\boldsymbol{y}_n,\Omega)} \to \chi_{\img(\boldsymbol{y},\Omega)} &\text{ in $L^1(\RN)$,}\\
		\label{eq:conv-inv}
		\overline{\boldsymbol{y}}_n^{-1}\to \overline{\boldsymbol{y}}^{-1} &\text{ in $L^1(\RN;\RN)$.}
	\end{align}
	Additionally, suppose that 
	\begin{equation}
		\label{eq:conv-adj}
		\text{$\adj \nabla \boldsymbol{y}_n \wk \adj \nabla \boldsymbol{y}$ in $L^1(\Omega;\RNN)$}
	\end{equation}
	and
	\begin{equation}
		\label{eq:equi-det}
		\text{$(\det \nabla \overline{\boldsymbol{y}}_n^{-1})_n$ is equi-integrable.}
	\end{equation}
	Then 
	\begin{equation}
		\label{eq:conv-nabla}
		\text{$ \nabla \overline{\boldsymbol{y}}_n^{-1}\wk \nabla \overline{\boldsymbol{y}}^{-1}$ in $L^1(\RN;\RNN)$.}
	\end{equation}
\end{lemma}
\begin{remark}[Equi-integrability]
	\label{rem:equi-det-inv}
Assumption \eqref{eq:equi-det} is satisfied when
\begin{equation*}
	\sup_{n\in\N} \int_\Omega \gamma(\det \nabla \boldsymbol{y}_n)\,\d\boldsymbol{x}<+\infty
\end{equation*}
for some Borel function $\gamma \colon (0,+\infty)\to [0,+\infty)$ satisfying the first condition in \eqref{eq:gamma}.
To see this, let $f\colon [0,+\infty)\to [0,+\infty)$ be defined by setting $f(s)\coloneqq s \gamma(1/s)$ if $s>0$ and $f(0)\coloneqq 0$. Then, the first condition in \eqref{eq:gamma}  superlinarity of $f$, namely
\begin{equation*}
	\lim_{s \to +\infty} \frac{f(s)}{s}=+\infty.
\end{equation*}
Using the change-of-variable formula, we compute
\begin{equation*}
	\begin{split}
		\int_{\RN} f(\det \nabla \overline{\boldsymbol{y}}_n^{-1})\,\d\boldsymbol{\xi}&=\int_{\img(\boldsymbol{y}_n,\Omega)} f(\det \nabla \boldsymbol{y}^{-1}_n)\,\d\boldsymbol{\xi}\\
		&=\int_{\img(\boldsymbol{y}_n,\Omega)} \gamma(\det \nabla \boldsymbol{y}_n)\circ \boldsymbol{y}_n^{-1}\det\nabla \boldsymbol{y}_n^{-1}\,\d\boldsymbol{\xi}=\int_\Omega \gamma(\det \nabla \boldsymbol{y}_n)\,\d\boldsymbol{x},
	\end{split}
\end{equation*}  
where the right-hand side is uniformly bounded with respect to $n\in\N$. Thus, the claim follows by the de la Vall\'{e}e Poussin criterion \cite[Theorem 2.24]{fonseca.leoni}.
 
\end{remark}
\begin{proof}
	Claims \eqref{eq:conv-img}--\eqref{eq:conv-inv} are established in \cite[Proposition 3.2(i)]{bresciani.friedrich.moracorral}. To prove \eqref{eq:conv-nabla}, we argue similarly to \cite[Theorem 6.3(c)]{barchiesi.henao.moracorral}. First, we  show that
	\begin{equation}
		\label{eq:equi-nabla}
		\text{$(\nabla \overline{\boldsymbol{y}}_n^{-1})_n$ is equi-integrable.}
	\end{equation}
	From \eqref{eq:conv-adj} and the Dunford-Pettis theorem, the sequence $(\adj \nabla \boldsymbol{y}_n)_n$ is equi-integrable.
	Let $\varepsilon>0$. Then,  there exists $\delta=\delta(\varepsilon)>0$ such that, for any measurable set $B\subset \Omega$,
	\begin{equation}
		\label{eq:equi-adj}
		\text{$\displaystyle \sup_{n\in\N} \int_B |\adj \nabla \boldsymbol{y}_n|\,\d\boldsymbol{x}<\varepsilon$ \quad whenever $\leb(B)<\delta$.}
	\end{equation} 
	Similarly, \eqref{eq:equi-det} yields $\omega=\omega(\delta)>0$ such that, for any measurable set $A\subset \RN$,
	\begin{equation}
		\label{eq:equi-det-inv}
		\sup_{n\in\N} \int_A \det \nabla \overline{\boldsymbol{y}}_n^{-1}\,\d\boldsymbol{\xi}<\delta \quad \text{whenever $\leb(A)<\omega.$}
	\end{equation} 
	Let $n\in\N$ and let $A\subset \RN$ be measurable with $\leb(A)<\omega$. Setting $B_n\coloneqq \boldsymbol{y}_n^{-1}(A\cap \img(\boldsymbol{y}_n,\Omega))$, using the area formula and \eqref{eq:equi-det-inv}, we find
	\begin{equation*}
		\leb(B_n)=\leb(\overline{\boldsymbol{y}}_n^{-1}(A))=\int_A \det \nabla \overline{\boldsymbol{y}}_n^{-1}\,\d\boldsymbol{\xi}<\delta
	\end{equation*} 
	and, in turn,
	\begin{equation*}
		\begin{split}
			\int_A |\nabla \overline{\boldsymbol{y}}_n^{-1}|\,\d \boldsymbol{\xi}&=\int_{A\cap \img(\boldsymbol{y}_n,\Omega)} |(\nabla \boldsymbol{y}_n)^{-1}\circ \boldsymbol{y}_n^{-1}|\,\d\boldsymbol{\xi}\\
			&=\int_{\boldsymbol{y}_n(B_n)} |(\nabla \boldsymbol{y}_n)^{-1}|\det \nabla \boldsymbol{y}_n\,\d\boldsymbol{x}=\int_{B_n} |\adj \nabla \boldsymbol{y}_n|\,\d\boldsymbol{x}<\varepsilon
		\end{split}
	\end{equation*}
	thanks to the change-of-variable formula and \eqref{eq:equi-adj}. This proves \eqref{eq:equi-det}. 
	
	From \eqref{eq:conv-img} and \eqref{eq:equi-nabla}, we find
	\begin{equation*}
		\lim_{n\to \infty} \int_{\RN \setminus \img(\boldsymbol{y},\Omega)} |\nabla \overline{\boldsymbol{y}}_n^{-1}|\,\d\boldsymbol{\xi}=\lim_{n\to \infty} \int_{\img(\boldsymbol{y}_n,\Omega) \setminus \img(\boldsymbol{y},\Omega)} |\nabla \overline{\boldsymbol{y}}_n^{-1}|\,\d\boldsymbol{\xi}=0. 
	\end{equation*}
	Hence, by the Dunford-Pettis theorem, the sequence $(\nabla \overline{\boldsymbol{y}}_n)_n$ is weakly precompact in $L^1(\RN;\RNN)$. Therefore, if we 
	show that 
	\begin{equation}
		\label{eq:nabla-cb}
		\text{$\nabla \overline{\boldsymbol{y}}_n^{-1}\wks \nabla \overline{\boldsymbol{y}}^{-1}$ in $C_{\rm b}^0(\RN;\RNN)'$,}
	\end{equation}
	then \eqref{eq:conv-nabla} follows in view of the Urysohn property. Let $\psi\in C_{\rm b}^0(\RN)$. By the change-of-variable formula, for all $n\in\N$,
	\begin{equation}
		\label{eq:psi}
		\int_{\RN} \psi\, \nabla\overline{\boldsymbol{y}}_n^{-1}\,\d\boldsymbol{\xi}=\int_{\img(\boldsymbol{y}_n,\Omega)} \psi\, (\nabla \boldsymbol{y}_n)^{-1}\circ \boldsymbol{y}_n^{-1}\,\d\boldsymbol{\xi}=\int_\Omega \psi\circ \boldsymbol{y}_n (\adj \nabla \boldsymbol{y}_n)\,\d\boldsymbol{x}
	\end{equation} 
	and, analogously,
	\begin{equation*}
		\int_{\RN} \psi\, \nabla\overline{\boldsymbol{y}}^{-1}\,\d\boldsymbol{\xi}=\int_\Omega \psi\circ \boldsymbol{y} (\adj \nabla \boldsymbol{y})\,\d\boldsymbol{x}.
	\end{equation*}
	Thus, passing to the limit, as $n\to \infty$, at the right-hand side of \eqref{eq:psi} with the aid of \cite[Proposition 2.61]{fonseca.leoni} having in mind the almost everywhere convergence in \eqref{eq:conv-ae-det} and \eqref{eq:conv-adj}, we obtain
	\begin{equation*}
		\lim_{n\to \infty} \int_{\RN} \psi\, \nabla\overline{\boldsymbol{y}}_n^{-1}\,\d\boldsymbol{\xi}=\int_{\RN} \psi\, \nabla\overline{\boldsymbol{y}}^{-1}\,\d\boldsymbol{\xi}.
	\end{equation*}
	This proves \eqref{eq:nabla-cb}.
\end{proof}

\subsection{Sobolev maps}
For the definition and the fundamental properties of the topological degree of continuous maps, we refer to \cite{fonseca.gangbo}. 
The next definition takes advantage of  Tiezte's extension theorem \cite[Theorem 1.15]{fonseca.gangbo} and the solely dependence of the degree on boundary values \cite[Theorem 2.4]{fonseca.gangbo}. Below, we denote by $\boldsymbol{y}^*$ the precise representative of a map $\boldsymbol{y}\in W^{1,p}(\Omega;\RN)$ defined as
\begin{equation*}
	\boldsymbol{y}^*(\boldsymbol{x})\coloneqq \limsup_{r\to 0^+} \dashint_{B(\boldsymbol{x},r)} \boldsymbol{y}(\boldsymbol{z})\,\d\boldsymbol{z}\quad \text{for all $\boldsymbol{x}\in\Omega$.}
\end{equation*}

\begin{definition}[Topological image]
	\label{def:topim}
Given $\boldsymbol{y}\in W^{1,p}(\Omega;\RN)$ and a domain $U\subset \subset \Omega$ such that $\boldsymbol{y}^*\restr{\partial U}\in C^0(\partial U;\RN)$, we define the topological image of $U$ under $\boldsymbol{y}$ as
\begin{equation*}
	\imt(\boldsymbol{y},U)\coloneqq \left\{ \boldsymbol{\xi}\in\RN \setminus {\boldsymbol{y}}^*(\partial U):\:\deg({\boldsymbol{y}}^*,U,\boldsymbol{\xi})\neq 0  \right\},
\end{equation*}
where  $\deg({\boldsymbol{y}}^*,U,\boldsymbol{\xi})$ denotes  the topological degree of any extension of ${\boldsymbol{y}}^*\restr{\partial U}$ to $\closure{U}$ at the point $\boldsymbol{\xi}$.  
Moreover, we define the topological image of  $\boldsymbol{y}$ as
\begin{equation*}
	\imt(\boldsymbol{y},\Omega)\coloneqq \bigcup \left\{ \imt(\boldsymbol{y},U):\:\text{$U\subset \subset \Omega$  domain with $\boldsymbol{y}^*\restr{\partial U}\in C^0(\partial U;\RN)$} \right\}.
\end{equation*}
\end{definition}

It is known that the topological degree defines a continuous function. As a consequence, the set $\imt(\boldsymbol{y},U)$ is  open and, in turn, so is $\imt(\boldsymbol{y},\Omega)$, with $\partial\, \imt(\boldsymbol{y},U)\subset \boldsymbol{y}^*(\partial U)$. Observe that $\deg(\boldsymbol{y}^*,U,\cdot)=0$ in the unique unbounded connected component of $\RN \setminus \boldsymbol{y}^*(\partial U)$, see, e.g., \cite[Remark 2.17(d)]{bresciani.friedrich.moracorral} for a proof. Thus,  ${\imt(\boldsymbol{y},U)}$ is bounded.  Because of the assumption $p>N-1$, if the boundary of $U$ is sufficiently regular and the trace of $\boldsymbol{y}$ on $\partial U$ belongs to $W^{1,p}(\partial U;\RN)$, then, by Morrey's embedding  \cite[Proposition A.10]{bresciani.friedrich.core},  this map admits a continuous representative. Using mollifications and the coarea formula, it can be shown that, for many domains $U$, the continuous representatives actually coincides with $\boldsymbol{y}^*\restr{\partial U}$.

Next, we present the invertibility condition (INV) due to {\sc M\"{u}ller \& Spector} \cite{mueller.spector}.

\begin{definition}[Invertibility condition]
\label{def:INV}	
Given $\boldsymbol{y}\in W^{1,p}(\Omega;\RN)$, we say that $\boldsymbol{y}$ satisfies {\rm (INV)} whenever, for all $\boldsymbol{a}\in\Omega$ and for almost all $r\in \left(0,\dist(\boldsymbol{a};\partial \Omega) \right)$, the following holds:
\begin{enumerate}[label=(\roman*)]
	\item $\boldsymbol{y}^*\restr{S(\boldsymbol{a},r)}\in C^0(S(\boldsymbol{a},r);\RN)$;
	\item $\boldsymbol{y}(\boldsymbol{x})\in\imt(\boldsymbol{y},B(\boldsymbol{a},r))$ for almost all $\boldsymbol{x}\in B(\boldsymbol{a},r)$;
	\item $\boldsymbol{y}(\boldsymbol{x})\notin\imt(\boldsymbol{y},B(\boldsymbol{a},r))$ for almost all $\boldsymbol{x}\in \Omega \setminus B(\boldsymbol{a},r)$.
\end{enumerate}
\end{definition}
By \cite[Lemma 3.4]{mueller.spector}, if $\det D\boldsymbol{y}(\boldsymbol{x})\neq 0$ for almost all $\boldsymbol{x}\in\Omega$ and $\boldsymbol{y}$ satisfies (INV), then $\boldsymbol{y}$ is almost everywhere injective. This fact, together with the approximate differentiability properties of Sobolev maps \cite[Theorem 2, p. 216]{gms.cc1}, gives the inclusion $\mathcal{Y}_p(\Omega)\subset \mathcal{Y}(\Omega)$.

To rigorously define cavities, we give the following definition.

\begin{definition}[Topological image of points and cavitation points] \label{def:cav}
Given $\boldsymbol{y}\in W^{1,p}(\Omega;\RN)$, for any $\boldsymbol{a}\in\Omega$, we define the topological image of $\boldsymbol{a}$ under $\boldsymbol{y}$ as
\begin{equation*}
	\imt(\boldsymbol{y},\boldsymbol{a})\coloneqq \bigcap \left\{ \closure{\imt(\boldsymbol{y},B(\boldsymbol{a},r))}:\:\text{$r\in \left(0,\dist(\boldsymbol{a};\partial\Omega)\right)$ with  $\boldsymbol{y}^*\restr{\partial U}\in C^0(S(\boldsymbol{a},r);\RN)$} \right\}.
\end{equation*} 	
Moreover, we define the set of cavitation points of $\boldsymbol{y}$ as
\begin{equation*}
	C_{\boldsymbol{y}}\coloneqq \left\{ \boldsymbol{a}\in\Omega:\:\leb(\imt(\boldsymbol{y},\boldsymbol{a}))>0  \right\}.
\end{equation*}
\end{definition} 

If $\boldsymbol{y}\in\mathcal{Y}_p(\Omega)$, then it has been proved in \cite[Lemma 7.3(i)]{mueller.spector} that $ \closure{\imt(\boldsymbol{y},B(\boldsymbol{a},r_1))}\subset\closure{\imt(\boldsymbol{y},B(\boldsymbol{a},r_2))}$ for any two radii $r_1$ and $r_2$ as above with $r_1<r_2$. As a consequence, $\imt(\boldsymbol{y},\boldsymbol{a})$ can be rewritten as a decreasing union of nonempty compact sets and, as such, it inherits the same two features.

Recall Definition \ref{def:surface-energy}. 
The following properties of deformations with finite surface energy have been essentially established in \cite{henao.moracorral.lusin}. We refer to \cite[Lemma 2.34 and Theorem 2.39]{bresciani.friedrich.moracorral} for the exact statement  below.

\begin{proposition}[Finite surface energy]
	\label{prop:finite-surface}
Let $\boldsymbol{y}\in\mathcal{Y}_p(\Omega)$ be such that $\mathcal{S}(\boldsymbol{y})<+\infty$. Then: 
\begin{enumerate}[label=(\roman*)]
	\item $ \imt(\boldsymbol{y},\Omega)\cong \img(\boldsymbol{y},\Omega)\cup \bigcup_{\boldsymbol{a}\in C_{\boldsymbol{y}}}  \imt(\boldsymbol{y},\boldsymbol{a})$, where  $\imt(\boldsymbol{y},\Omega) \cap \bigcup_{\boldsymbol{a}\in C_{\boldsymbol{y}}} \imt(\boldsymbol{y},\boldsymbol{a})\cong \emptyset$ and $\imt(\boldsymbol{y},\boldsymbol{a})\cap \imt(\boldsymbol{y},\boldsymbol{b})\simeq \emptyset$ for all distinct $\boldsymbol{a},\boldsymbol{b}\in C_{\boldsymbol{y}}$;
	 \item $\partial^*\img(\boldsymbol{y},\Omega)\cap \imt(\boldsymbol{y},\Omega)\simeq \bigcup_{\boldsymbol{a}\in C_{\boldsymbol{y}}} \partial^* \imt(\boldsymbol{y},\boldsymbol{a})$;
	 \item  $\mathcal{S}(\boldsymbol{y})=\sum_{\boldsymbol{a}\in C_{\boldsymbol{y}}} \per \left( \imt(\boldsymbol{y},\boldsymbol{a})\right)$, where the set $C_{\boldsymbol{y}}$ is countable.
\end{enumerate}	
\end{proposition}
\begin{remark}[Invertibility condition  and finite surface energy]\label{rem:inv-finite}
We stress that the validity of the previous proposition and, in particular the representation formula in (iii), crucially rely on the assumptions that $\boldsymbol{y}$ satisfies (INV) and $\mathcal{S}(\boldsymbol{y})<+\infty$. See \cite[pp. 51--54]{mueller.spector} for a counterexample and  also  \cite[pp. 629--631]{henao.moracorral.invertibility} for further comments in this direction. We mention that, Proposition \ref{prop:finite-surface} still holds true in the more general setting of deformations in the space $W^{1,N-1}(\Omega;\RN)\cap L^\infty(\Omega;\RN)$. We refer to \cite{henao.moracorral.lusin} for  details.
\end{remark}

For our purposes, the following definition of inverse is particularly convenient.  

\begin{definition}[Inverse deformation]
	\label{def:inverse-sobolev}
	Let $\boldsymbol{o}\notin\closure{\Omega}$ be arbitrarily fixed.
	Given $\boldsymbol{y}\in\mathcal{Y}_p(\Omega)$, we define the inverse deformation $\widehat{\boldsymbol{y}}^{-1}\colon \imt(\boldsymbol{y},\Omega)\to \RN$ of $\boldsymbol{y}$ by setting
	\begin{equation*}
		\widehat{\boldsymbol{y}}^{-1}(\boldsymbol{\xi})\coloneqq \begin{cases}
			\boldsymbol{y}^{-1}(\boldsymbol{\xi}) & \text{if $\boldsymbol{\xi}\in \img(\boldsymbol{y},\Omega)$,}\\
			\boldsymbol{o} & \text{if $\boldsymbol{\xi}\in \bigcup_{\boldsymbol{a}\in C_{\boldsymbol{y}}} \imt(\boldsymbol{y},\boldsymbol{a})$.}
		\end{cases} 
	\end{equation*} 
\end{definition}
In view of Proposition \ref{prop:finite-surface}(i), the previous definition specifies the value of $\widehat{\boldsymbol{y}}^{-1}$ almost everywhere in $\imt(\boldsymbol{y},\Omega)$. The choice of $\boldsymbol{o}$  makes sure that the jump of $\widehat{\boldsymbol{y}}^{-1}$ corresponds exactly to the boundary of the cavities. 

The next result on the regularity of inverse deformations as in Definition \ref{def:inverse-sobolev} is a simple corollary of \cite[Theorem 3.3]{henao.moracorral.regularity}. Further regularity results have been achieved in \cite[Theorem 2]{henao.moracorral.fracture} and \cite[Proposition 4]{henao.moracorral.xu}.

\begin{theorem}[Regularity of inverse deformations]
	\label{thm:regularity} 		
Let $\boldsymbol{y}\in\mathcal{Y}_p(\Omega)$ with $\mathcal{S}(\boldsymbol{y})<+\infty$. Then, $\widehat{\boldsymbol{y}}^{-1}\in SBV(\imt(\boldsymbol{y},\Omega);\RN)$ with distributional gradient
\begin{equation}
	\label{eq:SBV-gradient}
	D\widehat{\boldsymbol{y}}^{-1}=(\nabla\boldsymbol{y})^{-1}\circ \boldsymbol{y}^{-1} \leb \mres \img(\boldsymbol{y},\Omega)+\sum_{\boldsymbol{a}\in C_{\boldsymbol{y}}} (\boldsymbol{a}-\boldsymbol{o})\otimes \boldsymbol{\nu}_{\imt(\boldsymbol{y},\boldsymbol{a})}\haus \mres \partial^*\imt(\boldsymbol{y},\boldsymbol{a}). 
\end{equation}
In particular
\begin{equation}
	\label{eq:SBV-jump}
	J_{\widehat{\boldsymbol{y}}^{-1}}\simeq \bigcup_{\boldsymbol{a}\in C_{\boldsymbol{y}}} \partial^* \imt(\boldsymbol{y},\boldsymbol{a})
\end{equation}
and
\begin{equation}
	\label{eq:SBV-normal}
	\text{$\boldsymbol{\nu}_{\widehat{\boldsymbol{y}}_n^{-1}}\simeq -\boldsymbol{\nu}_{\imt(\boldsymbol{y},\boldsymbol{a})}$ on $\partial^*\imt(\boldsymbol{y},\boldsymbol{a})$ for all $\boldsymbol{a}\in C_{\boldsymbol{y}}$.}
\end{equation}
\end{theorem}
\begin{remark}[Normal to the geometric image]
Taking into account Proposition \ref{prop:finite-surface}(ii), it can be shown that $\boldsymbol{\nu}_{\widehat{\boldsymbol{y}}^{-1}}\simeq \boldsymbol{\nu}_{\img(\boldsymbol{y},\Omega)}$ on $\bigcup_{\boldsymbol{a}\in C_{\boldsymbol{y}}} \partial^*\imt(\boldsymbol{y},\boldsymbol{a})$. We refer to \cite[Theorem 4.8]{henao.moracorral.lusin} for related results.
\end{remark}
\begin{proof}
Define $\widetilde{\boldsymbol{y}}^{-1}\colon \imt(\boldsymbol{y},\Omega)\to \RN$ and $\boldsymbol{v}\colon \imt(\boldsymbol{y},\Omega)\to \RN$ by setting
\begin{equation*}
	\widetilde{\boldsymbol{y}}^{-1}(\boldsymbol{\xi})\coloneqq \begin{cases}
		\boldsymbol{y}^{-1}(\boldsymbol{\xi}) & \text{in $\img(\boldsymbol{y},\Omega)$,}\\
		\boldsymbol{a} & \text{in $\imt(\boldsymbol{y},\boldsymbol{a})$ for all $\boldsymbol{a}\in C_{\boldsymbol{y}}$},
	\end{cases} 
\end{equation*}
and
\begin{equation*}
	\boldsymbol{v}(\boldsymbol{\xi})\coloneqq \begin{cases}
		\boldsymbol{0} & \text{in $\img(\boldsymbol{y},\Omega)$,}\\
		\boldsymbol{o}-\boldsymbol{a} & \text{in $\imt(\boldsymbol{y},\boldsymbol{a})$ for all $\boldsymbol{a}\in C_{\boldsymbol{y}}$}.
	\end{cases}
\end{equation*}
These two functions are uniquely defined, up to representatives, in view of Definition \ref{def:inverse} and Proposition \ref{prop:finite-surface}(i).
Thanks to \cite[Theorem 3.3]{henao.moracorral.regularity}, we know that $\widetilde{\boldsymbol{y}}^{-1}\in W^{1,1}(\imt(\boldsymbol{y},\Omega);\RN)$ with $D\widetilde{\boldsymbol{y}}^{-1}= (\nabla \boldsymbol{y})^{-1}\circ \boldsymbol{y}^{-1}$ in $\img(\boldsymbol{y},\Omega)$ and $D\widetilde{\boldsymbol{y}}^{-1}= \boldsymbol{O} $ in $\bigcup_{\boldsymbol{a}\in C_{\boldsymbol{y}}} \imt(\boldsymbol{y},\boldsymbol{a})$.
From Proposition \ref{prop:finite-surface}, we realize that $\boldsymbol{v}\in SBV(\imt(\boldsymbol{y},\Omega);\RN)$ with distributional gradient
\begin{equation*}
	D\boldsymbol{v}=-\sum_{\boldsymbol{a}\in C_{\boldsymbol{y}}} (\boldsymbol{a}-\boldsymbol{o})\otimes \boldsymbol{\nu}_{\imt(\boldsymbol{y},\boldsymbol{a})}\haus\mres \partial^*\imt(\boldsymbol{y},\boldsymbol{a}).
\end{equation*}
Recalling Definition \ref{def:inverse-sobolev}, there holds $\widehat{\boldsymbol{y}}^{-1}=\widetilde{\boldsymbol{y}}^{-1}-\boldsymbol{v}$. 
Therefore, $\widehat{\boldsymbol{y}}^{-1}\in SBV(\imt(\boldsymbol{y},\Omega);\RN)$ and \eqref{eq:SBV-gradient} holds true. From this, \eqref{eq:SBV-jump}--\eqref{eq:SBV-normal} immediately follow.
\end{proof}

The next result provides a variant of Lemma \ref{lem:continuity-inverse} for  inverse deformations as in Definition \ref{def:inverse-sobolev}. 

\begin{proposition}[Continuity properties of the inverse]
	\label{prop:continuity-inverse}
Let $(\boldsymbol{y}_n)$ be a sequence in $\mathcal{Y}_p(\Omega)$ and $\boldsymbol{y}\in\mathcal{Y}_p(\Omega)$. Suppose that
\begin{equation}
	\label{eq:conv-sobolev}
	\text{$\boldsymbol{y}_n \wk \boldsymbol{y}$ in $W^{1,p}(\Omega;\RN)$,} \quad \text{$\det D \boldsymbol{y}_n \wk \det D \boldsymbol{y}$ in $L^1(\Omega)$.}
\end{equation} 
Also, let $V\subset \subset \imt(\boldsymbol{y},\Omega)$ be a bounded open set. Then, up to subsequences, $V\subset \subset \imt(\boldsymbol{y}_n,\Omega)$ for all $n\in \N$ and we have 
\begin{equation}
	\label{eq:conv-inverses}
	\text{$\widehat{\boldsymbol{y}}^{-1}_n \to  \widehat{\boldsymbol{y}}^{-1}$  in $L^1(V;\RN)$.}\\
\end{equation}
Additionally, suppose that
\begin{equation*}
	\text{$\adj D \boldsymbol{y}_n \wk \adj D\boldsymbol{y}$ in $L^1(\Omega;\RNN)$}
\end{equation*}
and
\begin{equation}
	\label{eq:equiintdet}
	\text{$(\det \nabla \overline{\boldsymbol{y}}_n^{-1})_n$ is equi-integrable.}
\end{equation}
Then
\begin{equation}
	\label{eq:conv-nabla-sobolev}
	\text{$\nabla \widehat{\boldsymbol{y}}_n^{-1} \wk \nabla \widehat{\boldsymbol{y}}^{-1}$ in $L^1(V;\RNN)$.}
\end{equation}
\end{proposition}
In  \eqref{eq:equiintdet}, we  refer to the extended inverses as in Definition \ref{def:inverse}.

\begin{remark}[Weak-* continuity of the inverse]
If we  know that $\boldsymbol{y}_n\wk \boldsymbol{y}$ in $\W^{1,p}(\Omega;\RN)$, $(\det D \boldsymbol{y}_n)_n$ is equi-integrable, and $\sup_{n\in\N} \mathcal{S}(\boldsymbol{y}_n)<+\infty$, then the weak continuity of the Jacobian determinant follows from \cite[Theorem]{henao.moracorral.invertibility}. Thus, the proposition yields \eqref{eq:conv-inverses}. If, additionally, $(\adj D \boldsymbol{y}_n)_n$ is bounded in $L^1(\Omega;\RNN)$, then, using Proposition \ref{prop:finite-surface}(iii), Theorem \ref{thm:regularity}, and the change-of-variable formula, we find 
\begin{equation*}
	\begin{split}
		|D\widehat{\boldsymbol{y}}_n^{-1}|(\imt(\boldsymbol{y}_n,\Omega))&\leq \|\nabla \boldsymbol{y}^{-1}_n\|_{L^1(\img(\boldsymbol{y}_n,\Omega))}+\left( \max \{ |\boldsymbol{x}|:\:\boldsymbol{x}+\Omega\} + |\boldsymbol{o}| \right) \sum_{\boldsymbol{a}\in C_{\boldsymbol{y}_n}} \per \left( \imt(\boldsymbol{y}_n,\boldsymbol{a})\right)\\
		&\leq \|\adj D\boldsymbol{y}_n\|_{L^1(\Omega;\RNN)}+ \left( \max \{ |\boldsymbol{x}|:\:\boldsymbol{x}+\Omega\} + |\boldsymbol{o}| \right) \mathcal{S}(\boldsymbol{y}_n),
	\end{split}
\end{equation*}
where the right-hand side is uniformly bounded for $n\in\N$. From this, we deduce that $(\widehat{\boldsymbol{y}}_n^{-1}\restr{V})_n$ is uniformly bounded in $BV(V;\RN)$. Therefore, we conclude that $\widehat{\boldsymbol{y}}_n^{-1} \wks \widehat{\boldsymbol{y}}^{-1}$ in $BV(V;\RN)$ thanks to   \eqref{eq:conv-inverses},  \cite[Proposition 3.13]{ambrosio.fusco.pallara}, and the Urysohn property.
\end{remark}

\begin{proof}
By \cite[Proposition 3.5(i)]{bresciani.friedrich.moracorral}, we have $V\subset \subset \imt(\boldsymbol{y}_n,\Omega)$ for all $n\in\N$ along a not relabeled subsequence. Hence, \eqref{eq:conv-inverses} and \eqref{eq:conv-nabla-sobolev} make sense. 
Given Definition \ref{def:inverse-sobolev}, 
 the remaining claims follow from Lemma \ref{lem:continuity-inverse} as $(\nabla \widehat{\boldsymbol{y}}_n^{-1})\restr{V}=(\nabla \overline{\boldsymbol{y}}_n^{-1})\restr{V}$ and $(\nabla \widehat{\boldsymbol{y}}^{-1})\restr{V}=(\nabla \overline{\boldsymbol{y}}^{-1})\restr{V}$.
\end{proof}

\section{Proof of the main {\MMM results}}
\label{sec:proof}

\MMM In this section, we present the proofs of our main results. First, we establish the existence result by applying \EEE 
 the direct method of the calculus of variations. Some of the arguments are classical and we repeat them just for completeness. The results in Section \ref{sec:prelim} will be now used for proving the lower semicontinuity of the surface term in \eqref{eq:E}.

\begin{proof}[Proof of Theorem \ref{thm:existence}]
The proof is subdivided into two steps. 

\emph{Step 1 (Compactness).} Set $m\coloneqq \inf\{ \mathcal{E}(\boldsymbol{y}):\:\boldsymbol{y}\in\mathcal{A}\}$. In view of the growth condition \eqref{eq:growth} in (W2) and ($\Phi$2),  we have $m\geq {\MMM \hat{c}}$. We assume $m<+\infty$ since otherwise the result is trivial.

Let $(\boldsymbol{y}_n)_n$ be a minimizing sequence for $\mathcal{E}$ in $\mathcal{A}$. Again from \eqref{eq:growth} and ($\Phi$2), we see that
 \begin{equation}
 	\label{eq:bounded-energy}
 	\begin{split}
 		{\MMM c} \sup_{n\in\N}	&  \|D\boldsymbol{y}_n\|^p_{L^p(\Omega;\RNN)} +  \sup_{n\in\N} \|\gamma(\det D \boldsymbol{y}_n)\|_{L^1(\Omega)}\\
 		&+ \MMM\hat{c} \EEE \sup_{n\in\N} \sum_{\boldsymbol{a}\in C_{\boldsymbol{y}_n}} \per \left( \imt(\boldsymbol{y}_n,\boldsymbol{a}) \right) \leq \sup_{n\in\N} \mathcal{E}(\boldsymbol{y}_n)<+\infty.
 	\end{split}
 \end{equation}
Thanks to the Poincar\'{e} inequality with trace term and the boundary condition given by $\boldsymbol{d}$, we obtain that  $(\boldsymbol{y}_n)_n$ is bounded in $W^{1,p}(\Omega;\RN)$. Also, recalling  the second condition in \eqref{eq:gamma},  the sequence $(\det D \boldsymbol{y})_n$ is equi-integrable in view of \eqref{eq:bounded-energy} and the De la Vall\'{e}e Poussin criterion. Thus, we find $\boldsymbol{y}\in W^{1,p}(\Omega;\RN)$ and $h\in L^1(\Omega)$ such that, up to subsequences, 
\begin{equation}
	\label{eq:h}
	\text{$\boldsymbol{y}_n\wk \boldsymbol{y}$ in $W^{1,p}(\Omega;\RN)$}, \qquad \text{$\det D \boldsymbol{y}_n \wk h$ in $L^1(\Omega)$.}
\end{equation} 
In particular, we see that $h(\boldsymbol{x})\geq 0$ for almost all $\boldsymbol{x}\in\Omega$. 
A standard contradiction argument based on the first condition in \eqref{eq:gamma}, shows that, actually, $h(\boldsymbol{x})>0$ for almost all $\boldsymbol{x}\in\Omega$. Namely, suppose that $A\coloneqq \{ \boldsymbol{x}\in \Omega:\:h(\boldsymbol{x})=0 \}$ satisfies $\leb(A)>0$. Then, \eqref{eq:h} gives $\det D \boldsymbol{y}_n \to 0$ in $L^1(A)$ and, up to subsequences, $\det D\boldsymbol{y}_n \to 0$ almost everywhere in $A$, so that \eqref{eq:gamma} entails $\gamma(\det D \boldsymbol{y}_n)\to +\infty$ almost everywhere in $A$. Thus, by applying Fatou's lemma, we find
\begin{equation*}
	+\infty=\int_A \liminf_{n\to \infty} \gamma(\det D \boldsymbol{y}_n)\,\d\boldsymbol{x}\leq \liminf_{n\to \infty} \int_A \gamma(\det D \boldsymbol{y}_n)\,\d\boldsymbol{x}.
\end{equation*}
{\MMM This contradicts \eqref{eq:bounded-energy} and proves that $\leb(A)=0$}.

By the weak continuity of Jacobian minors, we obtain
\begin{equation}
	\label{eq:a1}
	\text{$\adj D \boldsymbol{y}_n \wk \adj D \boldsymbol{y}$ in $L^{\frac{p}{N-1}}(\Omega;\RNN)$,}
\end{equation}
while, \eqref{eq:bounded-energy} and Proposition \ref{prop:finite-surface}(iii) trivially yield
\begin{equation*}
	\sup_{n\in\N} \mathcal{S}(\boldsymbol{y}_n)\leq \frac{1}{\MMM \hat{c}} \sup_{n\in\N} \mathcal{E}(\boldsymbol{y}_n)<+\infty
\end{equation*}
thanks to ($\Phi$2). 
Hence, in view of \cite[Theorem 3]{henao.moracorral.invertibility}, we conclude
\begin{equation}
	\label{eq:a2}
	h\cong \det D \boldsymbol{y}, \qquad \mathcal{S}(\boldsymbol{y})<+\infty.
\end{equation} 
By \cite[Lemma 3.3]{mueller.spector},  $\boldsymbol{y}$ satisfies  (INV) and, in turn, $\boldsymbol{y}\in\mathcal{A}$ by the weak continuity of the trace operator. 
 
\emph{Step 2 (Lower semicontinuity).}  The weak continuity of Jacobian minors together with \eqref{eq:h} and \eqref{eq:a2} gives 
\begin{equation*}
	\text{$\mathbf{M}(D\boldsymbol{y}_n)\wk\mathbf{M}(D\boldsymbol{y})$ in $L^1\left(\Omega;  \prod_{i=1}^{N-1} \R^{\binom{N}{i}\times \binom{N}{i}}\times (0,+\infty)\right)$,}
\end{equation*}
where $\mathbf{M}(\boldsymbol{F})$ collects all minors of $\boldsymbol{F}\in\RNN$. 
Recalling (W3), classical lower semicontinuity results yield
\begin{equation}
	\label{eq:w}
	\int_\Omega \MMM W(D\boldsymbol{y})\EEE\,\d\boldsymbol{x}\leq \liminf_{n\to \infty} \int_\Omega \MMM W(D\boldsymbol{y}_n)\EEE \,\d \boldsymbol{x}.
\end{equation}
For the surface term,
without loss of generality, we assume that its inferior limit, as $n\to \infty$, is attained as a limit.  Recall that $\widehat{\boldsymbol{y}}_n^{-1}\in SBV(\imt(\boldsymbol{y}_n,\Omega);\RN)$ and $\widehat{\boldsymbol{y}}^{-1}\in SBV(\imt(\boldsymbol{y},\Omega);\RN)$ in view of Theorem \ref{thm:regularity}.  Let $V\subset \subset \imt(\boldsymbol{y},\Omega)$ be a bounded open set. 
By Proposition \ref{prop:continuity-inverse}, up to subsequences, we have  $V\subset \subset \imt(\boldsymbol{y}_n,\Omega)$ for all $n\in\N$ and \eqref{eq:conv-inverses}. For a not relabeled subsequence, this actually yields 
\begin{equation}
	\label{eq:b1}
	\text{$\widehat{\boldsymbol{y}}_n^{-1}\restr{V}\to \widehat{\boldsymbol{y}}^{-1}\restr{V}$ in measure.}
\end{equation}
Also, note that 
\begin{equation}
	\label{eq:b2}
	\sup_{n\in\N} \|\widehat{\boldsymbol{y}}_n^{-1}\|_{L^\infty(V;\RN)}\leq \max\{ |\boldsymbol{x}|:\:\boldsymbol{x}\in\closure{\Omega}  \}<+\infty.
\end{equation}
From Remark \ref{rem:equi-det-inv}, we realize that  the sequence $(\det \nabla \overline{\boldsymbol{y}}_n^{-1})_n$ is equi-integrable.  By Proposition \ref{prop:continuity-inverse},  we obtain \eqref{eq:conv-nabla-sobolev}, so that
\begin{equation}
	\label{eq:b3}
	\text{$(\nabla \widehat{\boldsymbol{y}}_n^{-1}\restr{V})$ is equi-integrable}
\end{equation}
thanks to the Dunford-Pettis theorem. 
At this point,   by applying \cite[Theorem 3.6]{ambrosio} to the sequence $(\widehat{\boldsymbol{y}}_n^{-1}\restr{V})_n$ having in mind \eqref{eq:b1}--\eqref{eq:b3} and the assumptions ($\Phi$1)--($\Phi$3), we obtain
\begin{equation*}
	\int_{J_{\widehat{\boldsymbol{y}}^{-1}}\cap V} \phi(\boldsymbol{\nu}_{\widehat{\boldsymbol{y}}^{-1}})\,\d\haus \leq \liminf_{n\to \infty} 	\int_{J_{\widehat{\boldsymbol{y}}^{-1}_n}\cap V} \phi(\boldsymbol{\nu}_{\widehat{\boldsymbol{y}}^{-1}_n})\,\d\haus\leq \liminf_{n\to \infty} \int_{J_{\widehat{\boldsymbol{y}}^{-1}_n}} \phi(\boldsymbol{\nu}_{\widehat{\boldsymbol{y}}^{-1}_n})\,\d\haus.
\end{equation*}
Taking the supremum among all bounded open sets $V\subset \subset \imt(\boldsymbol{y},\Omega)$, we find
\begin{equation*}
	\int_{J_{\widehat{\boldsymbol{y}}^{-1}}} \phi(\boldsymbol{\nu}_{\widehat{\boldsymbol{y}}^{-1}})\,\d\haus \leq \liminf_{n\to \infty} \int_{J_{\widehat{\boldsymbol{y}}^{-1}_n}} \phi(\boldsymbol{\nu}_{\widehat{\boldsymbol{y}}^{-1}_n})\,\d\haus.
\end{equation*}
In view of Proposition \ref{prop:finite-surface}(i) and Theorem \ref{thm:regularity}, the previous inequality can be rewritten as
\begin{equation}
	\label{eq:f}
	\sum_{\boldsymbol{a}\in C_{\boldsymbol{y}}} \int_{\partial^*\imt(\boldsymbol{y},\boldsymbol{a})} \phi(\boldsymbol{\nu}_{\imt(\boldsymbol{y},\boldsymbol{a})})\,\d\haus \leq \liminf_{n\to \infty} \sum_{\boldsymbol{a}\in C_{\boldsymbol{y}_n}} \int_{\partial^*\imt(\boldsymbol{y}_n,\boldsymbol{a})} \phi(\boldsymbol{\nu}_{\imt(\boldsymbol{y}_n,\boldsymbol{a})})\,\d\haus. 
\end{equation}
Eventually, the combination of \eqref{eq:w} and \eqref{eq:f} gives
\begin{equation*}
	\mathcal{E}(\boldsymbol{y})\leq \liminf_{n\to \infty} \mathcal{E}(\boldsymbol{y}_n)=m,
\end{equation*}
which proves that $\boldsymbol{y}$ is a minimizer of $\mathcal{E}$ in $\mathcal{A}$.
\end{proof} 

\MMM Next, we derive our equilibrium equations by performing outer variations. Also here, some steps are standard. In the proof, we will employ an anisotropic version of the formula for the first variation of the perimeter to treat the surface term in \eqref{eq:E}. 

\begin{proof}[Proof of Theorem \ref{thm:condition}]
Let $\boldsymbol{y}$ be a minimizer of $\mathcal{E}$ in $\mathcal{A}$ and let $\boldsymbol{\psi}\in C^1_{\rm c}(\RN;\RN)$ satisfy $\boldsymbol{\psi}\circ \boldsymbol{d}=\boldsymbol{0}$ on $\Gamma$ in the sense of traces. For $t\in\R$ with $|t|\ll 1$, the function $\boldsymbol{h}_t\coloneqq \boldsymbol{id}+t\boldsymbol{\psi}$, where $\boldsymbol{id}$ denotes the identity map, is a diffeomorphism on $\RN$ with $\det D \boldsymbol{h}_t(\boldsymbol{\xi})>0$ for all $\boldsymbol{\xi}\in\RN$. 
By the chain rule (see, e.g., \cite[Lemma A.11]{bresciani.friedrich.moracorral}),  we have $\boldsymbol{h}_t\circ \boldsymbol{y}\in W^{1,p}(\Omega;\RN)$ with $D(\boldsymbol{h}_t \circ \boldsymbol{y})(\boldsymbol{x})=(D\boldsymbol{h}_t (\boldsymbol{y}(\boldsymbol{x})))D\boldsymbol{y}(\boldsymbol{x})$ and, in turn, $\det D (\boldsymbol{h}_t\circ \boldsymbol{y})(\boldsymbol{x})>0$,  for almost all $\boldsymbol{x}\in \Omega$. By arguing as in \cite[Proposition 5.1]{moracorral}, we see that $\boldsymbol{h}_t\circ \boldsymbol{y}$ satisfies (INV), $\mathcal{S}(\boldsymbol{y})<+\infty$, and the boundary condition $\boldsymbol{h}_t\circ \boldsymbol{y}=\boldsymbol{d}$ on $\Gamma$ holds in the sense of traces. Therefore, $\boldsymbol{h}_t\circ \boldsymbol{y}\in\mathcal{A}$ for $|t|\ll 1$, and, additionally, 
\begin{equation}
	\label{eq:hh}
	C_{\boldsymbol{h}_t\circ \boldsymbol{y}}=C_{\boldsymbol{y}}, \qquad \imt(\boldsymbol{h}_t\circ \boldsymbol{y},\boldsymbol{a})=\boldsymbol{h}_t(\imt(\boldsymbol{y},\boldsymbol{a})) \quad \text{for all $\boldsymbol{a}\in C_{\boldsymbol{y}}$.}
\end{equation}

Because of to the minimality of $\boldsymbol{y}$, we have 
\begin{equation}
	\label{eq:d1}
	\frac{\d}{\d t} \mathcal{E}(\boldsymbol{h}_t \circ \boldsymbol{y}) \Bigg |_{t=0}=0
\end{equation} 
as long as the derivative exists. For the elastic term, by arguing as in \cite[Theorem 2.3(i)]{ball.op}, we compute  its derivative as
\begin{equation}
	\label{eq:d2}
		\frac{\d}{\d t} \int_\Omega W(D\boldsymbol{y})\,\d\boldsymbol{x} \Bigg |_{t=0}=\int_\Omega \left( DW(D\boldsymbol{y})(D\boldsymbol{y}^\top)  \right):D\boldsymbol{\psi}(\boldsymbol{y})\,\d\boldsymbol{x}.
\end{equation}
For the derivative of the surface term, we employ an isotropic variant of the formula for the first variation of the anisotropic perimeter (see, e.g., \cite[Theorem 3.6]{bellettini.novaga.riey}): given a set $E\subset \RN$ with finite perimeter, the set $\boldsymbol{h}_t(E)$ has also finite perimeter for $|t|\ll 1$ and there holds
\begin{equation}
	\label{eq:fvap}
	\frac{\d}{\d t} \int_{\partial^* \boldsymbol{h}_t(E)} \phi(\boldsymbol{\nu}_{\boldsymbol{h}_t(E)})\,\d\haus \Bigg |_{t=0}= \int_{\partial^* E} \phi(\boldsymbol{\nu}_E)\,\div_{\phi}^{\partial^*E} \boldsymbol{\psi}\,\d\haus,
\end{equation}
where
\begin{equation*}
	\div_{\phi}^{\partial^* E} \boldsymbol{\psi}\coloneqq \div\,\boldsymbol{\psi}+\nabla \phi(\boldsymbol{\nu}_E)\cdot \frac{(D\boldsymbol{\psi})^\top \boldsymbol{\nu}_E}{\phi(\boldsymbol{\nu}_E)}.
\end{equation*}
For convenience, we write $E_{\boldsymbol{a}}\coloneqq \imt(\boldsymbol{y},\boldsymbol{a})$ for all  $\boldsymbol{a}\in C_{\boldsymbol{y}}$. Then, by applying \eqref{eq:hh} and \eqref{eq:fvap}, we find
\begin{equation}
	\label{eq:d3}
	\begin{split}
		\frac{\d}{\d t} \sum_{\boldsymbol{a}\in C_{\boldsymbol{y}}} \int_{\partial^* E_{\boldsymbol{a}}} \phi(\boldsymbol{\nu}_{E_{\boldsymbol{a}}})\,\d\haus \Bigg |_{t=0}&= \sum_{\boldsymbol{a}\in C_{\boldsymbol{y}}} \frac{\d}{\d t} \int_{\partial^* E_{\boldsymbol{a}}} \phi(\boldsymbol{\nu}_{E_{\boldsymbol{a}}})\,\d\haus \Bigg |_{t=0}\\
		&=\sum_{\boldsymbol{a}\in C_{\boldsymbol{y}}} \int_{\partial^* E_{\boldsymbol{a}}} \phi(\boldsymbol{\nu}_{E_{\boldsymbol{a}}})\, \div_\phi^{\partial^* E_{\boldsymbol{a}}}\boldsymbol{\psi}\,\d \haus. 
	\end{split}
\end{equation}
Here, in the first line, the exchange of derivative and summation is justified as there exists a constant $d>0$, depending only of $\boldsymbol{\psi}$, such that
\begin{equation*}
	\frac{1}{|t|} \left | \int_{\partial^* \boldsymbol{h}_t(E_{\boldsymbol{a}})} \phi(\boldsymbol{\nu}_{\boldsymbol{h}_t(E_{\boldsymbol{a}})})\,\d\haus - \int_{\partial^* E_{\boldsymbol{a}}} \phi(\boldsymbol{\nu}_{E_{\boldsymbol{a}}})\,\d\haus  \right |\leq d\,\per(E_{\boldsymbol{a}}) \quad \text{for all $\boldsymbol{a}\in C_{\boldsymbol{y}}$}
\end{equation*}
and
\begin{equation*}
	\sum_{\boldsymbol{a}\in C_{\boldsymbol{a}}} \per(E_{\boldsymbol{a}})\leq \frac{1}{\hat{c}} \sum_{\boldsymbol{a}\in C_{\boldsymbol{y}}} \int_{\partial^* E_{\boldsymbol{a}}} \phi(\boldsymbol{\nu}_{E_{\boldsymbol{a}}})\,\d\haus \leq \mathcal{E}(\boldsymbol{y})<+\infty
\end{equation*}
thanks to ($\Phi$2). Altogether, by combining \eqref{eq:d1}--\eqref{eq:d2} and \eqref{eq:d3}, we obtain
\begin{equation*}
	\int_\Omega \left( DW(D\boldsymbol{y})(D\boldsymbol{y})^\top \right):D\boldsymbol{\psi}(\boldsymbol{y})\,\d\boldsymbol{x}+ \sum_{\boldsymbol{a}\in C_{\boldsymbol{y}}} \int_{\partial^* E_{\boldsymbol{a}}} \phi(\boldsymbol{\nu}_{E_{\boldsymbol{a}}})\, \div_\phi^{\partial^* E_{\boldsymbol{a}}}\boldsymbol{\psi}\,\d \haus=0
\end{equation*}
which concludes the proof.
\end{proof}

\EEE

\section*{Acknowledgements}

The Author  acknowledges the support of the Alexander von Humboldt Foundation through the Humboldt Research Fellowship. He is member of GNAMPA (Gruppo Nazionale per l'Analisi Matematica, la Probabilit\'{a} e le loro Applicazioni) of INdAM (Istituto Nazionale di Alta \hbox{Matematica}).   \EEE

\vskip 10pt
{
	}

\end{document}